\documentclass[12pt]{amsart}
\usepackage[a4paper,margin=3cm]{geometry}
\usepackage{amsmath,amssymb,mathtools}
\usepackage{microtype}
\usepackage[hidelinks]{hyperref}
\newcommand{\E}{\mathbb E}
\newcommand{\N}{\mathbb N}
\newcommand{\Z}{\mathbb Z}
\newcommand{\C}{\mathbb C}
\newcommand{\1}{\mathbf 1}
\newcommand{\e}{\mathrm e}

\newtheorem{theorem}{Theorem}[section]
\newtheorem{proposition}[theorem]{Proposition}
\newtheorem{lemma}[theorem]{Lemma}
\newtheorem{corollary}[theorem]{Corollary}
\newtheorem{remark}[theorem]{Remark}

\title{Polynomial configurations and pointwise averages
along Piatetski--Shapiro primes}
\author[H. Mousavi]{Hamed Mousavi}
\address{Department of Mathematics, University of Bristol,
Bristol BS8 1QU, UK}
 \email{gj23799@bristol.ac.uk}
\date{}

\begin{document}

\begin{abstract}
For every \(c>1\) sufficiently close to \(1\), we prove that there is \(\kappa>0\) such that every sufficiently large \(A\subset[N]\) of density at least \((\log N)^{-\kappa}\) contains
\[
x,\ x+(p-1),\ x+(p-1)^2,\ldots,\ x+(p-1)^k,
\qquad p=\lfloor n^c\rfloor\in\mathbb P.
\]
We also prove pointwise almost-everywhere convergence of the corresponding multiple ergodic averages.
%, with the same limit as the averages along the ordinary primes.
\end{abstract}

\maketitle
\pagestyle{plain}

\section{Introduction and transference}

The polynomial Szemer\'edi theorem of Bergelson and Leibman
\cite{BergelsonLeibman} generalizes S\'ark\"ozy's
square-difference theorem \cite{SarkozySquares} to polynomial
configurations in sets of positive upper density.
Wooley and Ziegler restricted the common parameter to shifted primes
\cite{WooleyZiegler}. Building on these results, Krause, the author,
Tao and Ter\"av\"ainen obtained the quantitative density bound
$(\log N)^{-\kappa}$ for distinct-degree polynomial families
\cite[Theorem~1.2]{KMTT}.

Szemer\'edi- and Green--Tao-type problems for arithmetically
defined subsets of the primes have attracted particular
interest when these sets have power-saving zero relative
density. Examples include Piatetski--Shapiro primes
\cite{MirekPS,LiPanPS,ShaoSun} and primes with restricted
digits, studied by Maynard \cite{MaynardDigits}, for which
Burgin obtained related recurrence results \cite{BurginDigits}.
% We ask whether polynomial
% Szemer\'edi-type configurations can still be forced when the prime
% parameter is restricted to a set this sparse.

For $c>1$, put
\[
 \mathcal S_c:=\{\lfloor n^c\rfloor:n\in\N\},\qquad
 \mathbb P_c:=\mathcal S_c\cap\mathbb P.
\]
When $c$ is sufficiently close to $1$,
\[
 \frac{|\mathbb P_c\cap[N]|}{|\mathbb P\cap[N]|}
 \sim N^{-(1-1/c)}.
\]
This follows from the Piatetski--Shapiro prime number theorem
\cite{PiatetskiShapiro}, which gives
$|\mathbb P_c\cap[N]|\sim N^{1/c}/\log N$, together with the
ordinary prime number theorem
$|\mathbb P\cap[N]|\sim N/\log N$.
Thus $\mathbb P_c$ has zero relative density among the ordinary
primes, with a power saving. Our first theorem shows that, for distinct-degree
polynomial families, restricting the gap size parameter to $\mathbb P_c-1$
still permits a density threshold $(\log N)^{-\kappa}$ for the integer
set containing the configuration. Up to a change in $\kappa>0$, this
matches the form of Shao and Wang's bound for unrestricted integer
parameters \cite[Theorem~1.1]{ShaoWang} and the bound of Krause, the
author, Tao and Ter\"av\"ainen for shifted-prime parameters
\cite[Theorem~1.2]{KMTT}.

Fix throughout
\begin{equation}\label{e:polynomials}
 \begin{gathered}
 P_0=0,\qquad P_1,\ldots,P_\ell\in\Z[y],\qquad P_i(0)=0,\qquad
 1\leq\deg P_1<\cdots<\deg P_\ell=d.
 \end{gathered}
\end{equation}
We write $[N]=\{1,\ldots,N\}$,
$\E_{n\leq N}=N^{-1}\sum_{n\leq N}$, and
$\mathbf P=(P_1,\ldots,P_\ell)$. Unrestricted sums over an integer
variable are over $\Z$. 
% The proof uses the distinct-degree generalized von Neumann theorem of
% Ter\"av\"ainen. 

% We impose the further requirement that the prime belong to the sparse
% Piatetski--Shapiro sequence.  For $c>1$, write
% \[
%  \mathcal S_c:=\{\lfloor n^c\rfloor:n\in\N\},\qquad
%  \mathbb P_c:=\mathcal S_c\cap\mathbb P.
% \]
% Shao and Sun recently proved that the normalized Piatetski--Shapiro-prime
% weight is asymptotically indistinguishable from the von Mangoldt function
% against arbitrary fixed-step nilsequences \cite[Theorem~2.1]{ShaoSun}.
% Among other consequences, they obtained asymptotic formulas for
% finite-complexity systems of linear equations in $\mathbb P_c$.  We use
% the polynomial-phase consequence of their theorem to transfer the
% distinct-degree results along the shifted primes to $\mathbb P_c$.

% Fix throughout
% \begin{equation}\label{e:polynomials}
%  P_0=0,\qquad P_1,\ldots,P_\ell\in\Z[y],\qquad
%  P_i(0)=0,\qquad
%  1\leq\deg P_1<\cdots<\deg P_\ell=d.
% \end{equation}
% We write $[N]=\{1,\ldots,N\}$,
% $\E_{n\leq N}=N^{-1}\sum_{n\leq N}$, and
% $\mathbf P=(P_1,\ldots,P_\ell)$.  Unrestricted sums over an integer
% variable are over $\Z$.

\begin{theorem}
\label{t:szemeredi}
For every family of polynomials satisfying~\eqref{e:polynomials},
there are $c_*>1$ and $\kappa>0$ with the following property.
For every $1<c<c_*$, there exists $N_0=N_0(\mathbf P,c)$ such
that, for every integer $N\geq N_0$, every $A\subset[N]$
satisfying
\[
 |A|\geq\frac{N}{(\log N)^\kappa}
\]
contains a configuration
\begin{equation}\label{e:pattern}
 x,\quad x+P_1(p-1),\quad\ldots,\quad x+P_\ell(p-1)
\end{equation}
for some $p\in\mathbb P_c$, with the displayed terms
pairwise distinct.
\end{theorem}

The connection between polynomial configurations and multiple recurrence
suggests a corresponding pointwise question: do simultaneous polynomial
returns have a limiting frequency along almost every orbit when the
prime parameter is restricted to a set of power-saving zero relative
density?

For ordinary primes, pointwise convergence of single linear averages
goes back to Bourgain and Wierdl \cite{BourgainPrimes,WierdlPrimes}.
Wooley and Ziegler proved norm convergence for multiple polynomial
averages along primes \cite{WooleyZiegler}. Pointwise convergence was
subsequently obtained by Krause, the author, Tao and Ter\"av\"ainen for
the bilinear family $(n,P(n))$, $\deg P\geq2$
\cite{KMTTPointwise}, and by Mirek, Wan and Wright for general
distinct-degree families \cite{MirekWanWright}.

For thin prime sets, Mirek established single pointwise ergodic
theorems, including along Piatetski--Shapiro primes
\cite{MirekThinPrimes}. Fornal and Krause treated bilinear linear
averages along primes and their Piatetski--Shapiro restrictions
\cite{FornalKrause}. The author proved a distinct-degree pointwise
theorem along the Piatetski--Shapiro integers in
\cite[Theorem~1.2]{mousavisolo}. Our next theorem shows that, for
$c>1$ sufficiently close to $1$, restricting to the polynomially
sparser set $\mathbb P_c$ preserves both pointwise convergence
and the limit of the ordinary-prime averages.  

 Let $(X,\nu,T)$ be an invertible
probability measure-preserving system and let
$f_1,\ldots,f_\ell\in L^\infty(X)$.  Put
\begin{equation}\label{e:PS-average}
 A_N^{\mathbb P_c}(\mathbf f)(x)
 :=\frac1{\pi_c(N)}
 \sum_{\substack{p\leq N\\p\in\mathbb P_c}}
 \prod_{i=1}^{\ell}f_i\bigl(T^{P_i(p-1)}x\bigr),
 \qquad \pi_c(N):=|\mathbb P_c\cap[N]|.
\end{equation}

\begin{theorem}[Pointwise averages along Piatetski--Shapiro primes]
\label{t:pointwise}
After decreasing $c_*>1$ if necessary, let $1<c<c_*$.  Then the
averages~\eqref{e:PS-average} converge almost everywhere, and their limit
agrees almost everywhere with the limit of
\begin{equation}\label{e:prime-average}
 \frac1{\pi(N)}\sum_{p\leq N}
 \prod_{i=1}^{\ell}f_i\bigl(T^{P_i(p-1)}x\bigr).
\end{equation}
Here $\pi(N)=|\mathbb P\cap[N]|$.
\end{theorem}

The comparison underlying this theorem also implies lacunary
$r$-variation bounds for the weighted averages when $r>2$;
see Corollary~\ref{c:long-variation}.

\begin{remark}
We expect Theorem~\ref{t:pointwise} to extend to
$f_i\in L^{p_i}(X)$ with $p_i>1$ and
$\sum_{i=1}^{\ell}p_i^{-1}\leq1$, allowing $c_*$ to depend on
the exponents. Such extensions lie beyond the scope of this paper.
The $L^1$ endpoint is more delicate, even in the single-function
case: LaVictoire proved that pointwise convergence can fail for
linear averages along the ordinary primes \cite{LaVictoireBad}.
Mirek raised the corresponding endpoint question for thin
prime sets \cite[discussion following Theorem~1.14]{MirekThinPrimes}.
\end{remark}

Both arguments begin by approximating the normalized
Piatetski--Shapiro-prime weight by the ordinary prime weight.
Writing $\gamma=1/c$ and letting $\Lambda$ denote the von Mangoldt
function, put
\begin{align}\label{e:weights}
    \Lambda_\gamma(n)
 :=\gamma^{-1}n^{1-\gamma}\Lambda(n)\1_{\mathcal S_c}(n),
 \qquad D_\gamma:=\Lambda_\gamma-\Lambda.
\end{align}
The normalization compensates for the sparsity of $\mathbb P_c$.
Centering in this way reduces both problems to showing that the
error weight $D_\gamma$ contributes negligibly to the relevant
polynomial averages.

A customary approach, based on PET induction and generalized von
Neumann estimates, controls this contribution by a Gowers norm of
the error weight. Here it is important to retain a power saving
despite the growing normalization factor $n^{1-\gamma}$. We instead
exploit the distinct degrees: Ter\"av\"ainen's theorem controls the
error directly through its correlations with ordinary polynomial
phases, avoiding the need for a quantitative Gowers-uniformity
estimate for $D_\gamma$.

Write $\e(t)=e^{2\pi it}$ and define
\begin{equation}\label{e:u-norm}
 \|w\|_{u^s[M]}
 :=\sup_{\substack{\phi\in\mathbb R[y]\\\deg\phi<s}}
 \left|\E_{m\leq M}w(m)\e(-\phi(m))\right|.
\end{equation}
A function is called $1$-bounded if its absolute value is at most $1$. 

\begin{theorem}[Ter\"av\"ainen,
 {\cite[Theorem~4.1]{Teravainen}}]\label{t:teravainen}
Let $R_1,\ldots,R_\ell\in\Z[y]$ satisfy
$\deg R_1<\cdots<\deg R_\ell=d$, and put
$\mathbf R=(R_1,\ldots,R_\ell)$. Let $C_0\geq1$.
There is an integer $1\leq K\ll_d1$ such that, whenever
$F_0,\ldots,F_\ell:\Z\to\C$ are $1$-bounded and supported on
$[-C_0M^d,C_0M^d]$, and $\theta:[M]\to\C$ is $1$-bounded,
\begin{equation}\label{e:teravainen}
 \begin{aligned}
 &\left|\frac1{M^{d+1}}\sum_{x\in\Z}\sum_{m\leq M}
 \theta(m)F_0(x)\prod_{i=1}^{\ell}F_i(x+R_i(m))\right|\ll_{C_0,\mathbf R}
 \bigl(M^{-1}+\|\theta\|_{u^{d+1}[M]}\bigr)^{1/K}.
 \end{aligned}
\end{equation}
\end{theorem}

By homogeneity, the same average with $\theta$ replaced by a weight
$g$ satisfying $|g|\leq B$, $B>0$, is bounded, up to the same implied
constant, by
\begin{equation}\label{e:scaled-vn}
 B\left(M^{-1}+B^{-1}\|g\|_{u^{d+1}[M]}\right)^{1/K}.
\end{equation}
Shao and Sun's estimate \cite[Theorem~2.1]{ShaoSun} supplies a power
saving for the correlations of $D_\gamma$ with polynomial phases,
uniformly in their real coefficients. 
% Taking $c>1$ sufficiently
% close to $1$ therefore gives a power-saving comparison between
% the sparse-prime and ordinary-prime averages.

For Theorem~\ref{t:szemeredi}, we insert this comparison into the
density increment argument of \cite{KMTT}. Its PET induction comparison
estimates and density increments accommodate the customary
$(C,Q)$-coefficient condition, which records the coefficient growth
caused by restricting the original set to arithmetic progressions. In our iteration,
the resulting polynomials have the special form
\begin{equation}\label{e:rescaled-polynomial}
 P_{i,Q}(y):=\frac{P_i(Qy)}Q.
\end{equation}
Applying Theorem~\ref{t:teravainen} directly to this changing family
would introduce constants depending on $Q$. Instead, the change of
variables $n=Qm$, $z=Qx$, with the functions supported on multiples
of $Q$, restores the original fixed polynomials $P_i$.
The resulting loss is only a power of $Q$. Since $Q=N^{o(1)}$
throughout the iteration, Shao--Sun's power saving absorbs this
loss. Thus no new coefficient-uniform version of
Theorem~\ref{t:teravainen} is needed. The comparison allows the
absence of a sparse-prime configuration to force a density increment;
repeated increments then give the required contradiction.

For Theorem~\ref{t:pointwise}, transference gives a power-saving
$L^1$ estimate for the difference between the two weighted ergodic
averages.  The ordinary-prime averages
converge by \cite[Theorem~1.1]{MirekWanWright}, so the sparse-prime
averages have the same limit. Finally, removing higher prime powers
and applying partial summation gives the unweighted version of the ergodic theorem.

Section~2 proves the resulting comparison with ordinary prime counts.
Section~3 inserts that comparison into the density increment of
\cite{KMTT}.  Section~4 proves Theorem~\ref{t:pointwise}.

\subsection*{Acknowledgements}

We thank Yu-Chen Sun and Anastasios Fragkos for helpful discussions.
ChatGPT was used as an auxiliary editorial tool for proofreading
and preliminary literature searches, and to generate an accompanying
Lean development~\cite{MousaviPSFormalization}. This development formalizes
the deductions within the paper conditional on eight explicitly
documented external analytic and arithmetic input interfaces.
The external results and the remaining source-to-interface
translations are not themselves proved in Lean.
All mathematical statements, references, arguments, and final wording
in the manuscript were checked and approved by the author, who takes
full responsibility for its contents. The Lean development is
publicly available but has not yet undergone independent human
expert review.

\section{The Piatetski--Shapiro-prime comparison}

Recall the notation of \eqref{e:weights} and extend the weights
by zero to non-positive integers. Contributions from higher prime
powers will be removed by elementary estimates.
We first record the arithmetic input. 

\begin{theorem}[Shao--Sun polynomial-phase estimate]
\label{t:shao-sun}
For every $r\geq1$ there are $\varepsilon_r,\sigma_r>0$ such that, if
$1-\varepsilon_r<\gamma<1$, then
\begin{equation}\label{e:shao-sun}
 \sup_{\deg\phi\leq r}
 \left|\sum_{n\leq X}D_\gamma(n)\e(\phi(n))\right|
 \ll_{r,\gamma}X^{1-\sigma_r}.
\end{equation}
\end{theorem}

The case $r=1$ goes back to Balog--Friedlander
\cite{BalogFriedlander}, while polynomial phases of general
degree were studied by Akbal--G\"ulo\u{g}lu
\cite{AkbalGuloglu}.
We obtain~\eqref{e:shao-sun} from the Shao--Sun nilsequence
estimate \cite[Theorem~2.1]{ShaoSun}.

\begin{proof}
Take the additive group $G=\mathbb R$, the lattice $\Gamma=\mathbb Z$,
and the filtration
\[
 G_0=\cdots=G_r=\mathbb R,\qquad G_j=\{0\}\quad(j>r).
\]
With the standard basis, this is a filtered nilmanifold of degree $r$,
dimension $1$, and complexity $O_r(1).$ The function
$F(t+\mathbb Z):=\e(t)$ has Lipschitz norm $O(1)$.

Every real polynomial $\phi$ of degree at most $r$ admits an expansion
\[
 \phi(n)=\sum_{j=0}^{r}a_j\binom nj,\qquad a_j\in\mathbb R.
\]
Thus $g_\phi(n):=\phi(n)$ is a polynomial sequence adapted to this
filtration, in the sense of \cite[Definition~4.5]{ShaoSun}.
Choose a fixed $\delta_r>0$, depending only on $r$, such that the
nilmanifold complexity and $\|F\|_{\mathrm{Lip}}$ are at most
$\delta_r^{-1}$.

Applying \cite[Theorem~2.1]{ShaoSun} with degree $r$, dimension $1$,
and these fixed data gives $\varepsilon_r,\sigma_r>0$ such that,
whenever $1-\varepsilon_r<\gamma<1$,
\[
 \left|\sum_{n\leq X}
 \bigl(\Lambda_\gamma(n)-\Lambda(n)\bigr)
 F(g_\phi(n)+\Gamma)\right|
 \ll_{r,\gamma}X^{1-\sigma_r}.
\]
Here $F(g_\phi(n)+\Gamma)=\e(\phi(n))$. The estimate is uniform in
$g_\phi$: the cited theorem imposes no coefficient bound on the
polynomial sequence, while all other data above depend only on $r$.
Taking the supremum over $\deg\phi\leq r$ proves
\eqref{e:shao-sun}.
\end{proof}
The next lemma restricts the estimate to the progressions arising in the
density increment.

\begin{lemma}\label{l:progression-phase}
For every integer $r\geq1$, there are $\varepsilon_r,\sigma_r>0$
such that, whenever $1-\varepsilon_r<\gamma<1$ and
$Q,M\geq1$ are integers, one has
\begin{equation}\label{e:progression-phase}
 \sup_{\substack{\psi\in\mathbb R[y]\\\deg\psi\leq r}}
 \left|\sum_{y\leq M}D_\gamma(Qy+1)\e(\psi(y))\right|
 \ll_{r,\gamma}(QM)^{1-\sigma_r}.
\end{equation}
\end{lemma}

\begin{proof}
Write $L\coloneqq QM$. The change of variables $n=Qy+1$ implies that  
\[
 \sum_{y\leq M}D_\gamma(Qy+1)\e(\psi(y))
 =
 \sum_{\substack{2\leq n\leq L+1\\n\equiv1\pmod Q}}
 D_\gamma(n)\e\!\left(\psi\!\left(\frac{n-1}{Q}\right)\right).
\]
By orthogonality,
\[
 \1_{\{n\equiv1\; (\mathrm{mod}\;Q)\}}
 =\frac1Q\sum_{a\bmod Q}\e\left(\frac{a(n-1)}Q\right).
\]
For each $a$, the phase
\[
 \psi\left(\frac{n-1}Q\right)+\frac{a(n-1)}Q
\]
is a real polynomial in $n$ of degree at most $r$.
Since $D_\gamma(1)=0$, we may include $n=1$ in each inner sum.
Theorem~\ref{t:shao-sun}, applied with $X=L+1$, gives
\[
 \sup_{a\bmod Q}
 \left|\sum_{n=1}^{L+1}D_\gamma(n)
 \e\!\left(\psi\!\left(\frac{n-1}{Q}\right)
 +\frac{a(n-1)}Q\right)\right|
 \ll_{r,\gamma}L^{1-\sigma_r}.
\]
Averaging over the $Q$ additive characters and using the
uniformity in $\psi$ proves the lemma.
\end{proof}

For the fixed family $\mathbf P$ in~\eqref{e:polynomials}, recall
the rescaled polynomials $P_{i,Q}$ from~\eqref{e:rescaled-polynomial}.
For a weight $w:\N\to\C$ and finitely supported functions
$\mathbf f=(f_0,\ldots,f_\ell)$ on $\Z$, put
\begin{equation}\label{e:counting-form}
 \mathcal T_{N,M,Q}(w;\mathbf f)
 :=\frac1N\sum_{x\in\Z}\E_{y\leq M}
 w(Qy+1)f_0(x)\prod_{i=1}^{\ell}
 f_i\bigl(x+P_{i,Q}(y)\bigr).
\end{equation}

\begin{proposition}
\label{p:comparison}
Fix the family~\eqref{e:polynomials}, $C_0\geq1$, and $0<a<1$.
There are constants $c_*=c_*(d)>1$ and $\xi=\xi(d)>0$
with the following property.
Let $1<c<c_*$, put $\gamma=1/c$, and let $N,Q$ be positive
integers with $Q\leq\exp((\log N)^a)$. Set
\[
 M:=\left\lfloor\left(\frac{N}{Q^{d-1}}\right)^{1/d}\right\rfloor.
\]
If $M$ is sufficiently large in terms of $c,\mathbf P,C_0,a$
and $f_0,\ldots,f_\ell$ are $1$-bounded and supported on
$[-C_0N,C_0N]$, then
\begin{equation}\label{e:comparison}
 \left|\mathcal T_{N,M,Q}(D_\gamma;\mathbf f)\right|
 \ll_{c,\mathbf P,C_0}N^{-\xi}.
\end{equation}
\end{proposition}

\begin{proof}
Put $L=QM$.  Since $Q=N^{o(1)}$,
\begin{equation}\label{e:lifting-scales}
 L^d\asymp QN,\qquad
 \frac{L^{d+1}}{NM}\asymp Q^2.
\end{equation}
For $0\leq i\leq\ell$, define
\[
 F_i(Qx)=f_i(x),\qquad F_i(z)=0\quad\text{if }Q\nmid z.
\]
These functions are $1$-bounded and supported on
$[-C_1L^d,C_1L^d]$ for a constant $C_1=C_1(C_0,\mathbf P)$.
Furthermore,
\begin{equation}\label{e:lift-identity}
 F_i\bigl(Qx+P_i(Qy)\bigr)
 =f_i\bigl(x+P_{i,Q}(y)\bigr).
\end{equation}
For $n\leq L$, set
\[
 \Theta(n):=B^{-1}D_\gamma(n+1)\1_{\{Q\mid n\}},
 \qquad B:=C_2L^{1-\gamma}\log(2L),
\]
where $C_2=C_2(\gamma)$ is chosen so that $|\Theta|\leq1$.  By
Lemma~\ref{l:progression-phase},
\begin{align}\label{e:theta-u}
 \|\Theta\|_{u^{d+1}[L]}
 &=\frac1{BL}
 \sup_{\deg\phi\leq d}
 \left|\sum_{y\leq M}D_\gamma(Qy+1)\e(-\phi(Qy))\right|\ll_{c,d}B^{-1}L^{-\sigma_d}.
% \nonumber
\end{align}
Note that every polynomial phase in $y$ occurs as $\phi(Qy)$ for a
real polynomial $\phi$, so no restriction of the supremum is hidden in
this identity.
Using~\eqref{e:lift-identity}, \eqref{e:theta-u}, and the definitions of $F_i$ and $\Theta$,
followed by Theorem~\ref{t:teravainen}, we obtain
\begin{align*}
 \left|\mathcal T_{N,M,Q}(D_\gamma;\mathbf f)\right|
 &=\frac{B}{NM}
 \left|\sum_{z\in\Z}\sum_{n=1}^{L}
 \Theta(n)F_0(z)\prod_{i=1}^{\ell}F_i(z+P_i(n))\right|\\
 &\ll_{c,\mathbf P,C_0}
 B\frac{L^{d+1}}{NM}
 \left(L^{-1}+B^{-1}L^{-\sigma_d}\right)^{1/K}\\
 &\ll_{c,\mathbf P,C_0}
 BQ^2\left(L^{-1}+B^{-1}L^{-\sigma_d}\right)^{1/K}.
\end{align*}
To make the final saving explicit, choose $\gamma$ so close to $1$ that
\[
 1-\gamma<\frac{\min(1,\sigma_d)}{8K}.
\]
For sufficiently large $L$ we then have $B\ll_c L^{2(1-\gamma)}$.
Splitting the last parenthesis into its two terms shows that 
\[
 \left|\mathcal T_{N,M,Q}(D_\gamma;\mathbf f)\right| \ll_{c,\mathbf P,C_0}Q^2\left(
 L^{2(1-\gamma)-1/K}
 +L^{2(1-\gamma)(1-1/K)-\sigma_d/K}
 \right).
\]
Both powers of $L$ are negative.  Since $Q=L^{o(1)}$ and
$L=N^{1/d+o(1)}$, this is $O_{c,\mathbf P,C_0}(N^{-\xi})$ for some
$\xi>0$.
\end{proof}

% \begin{remark}
% The comparison uses Theorem~\ref{t:teravainen} only for the original
% fixed family $P_1,\ldots,P_\ell$.  The factor $Q^2$ in the proof is the
% cost of restricting both variables to multiples of $Q$; it is harmless
% because all moduli in the density increment are $N^{o(1)}$.
% \end{remark}

\section{Quantitative recurrence}

We now insert Proposition~\ref{p:comparison} into the distinct-degree
density increment of \cite{KMTT}.  We state the exact interface that we
use.

For an integer $s\geq2$, write
$\mathcal L(s):=\operatorname{lcm}(1,\ldots,s)$, and let $P^+(Q)$
be the largest prime factor of $Q$, with $P^+(1)=1$.
For $C\geq1$ and an integer $Q\geq1$, an integer polynomial
$R(y)=a_ey^e+\cdots+a_1y$ of degree $e\geq1$ has
$(C,Q)$-coefficients if there exists $a'_e\in\Z$ such that
\[
 a_e=a'_eQ^{e-1},\qquad 0<|a'_e|\leq C,\qquad
 |a_j|\leq C|a_e|\quad(1\leq j<e).
\]
For the fixed family~\eqref{e:polynomials}, choose $C\geq1$
bounding the absolute values of all coefficients of
$P_1,\ldots,P_\ell$. Then $P_{i,Q}$ has $(C,Q)$-coefficients
for every integer $Q\geq1$ and $1\leq i\leq\ell$.
Indeed, writing $d_i=\deg P_i$ and
$P_i(y)=\sum_{j=1}^{d_i}a_{i,j}y^j$, we have
\begin{equation}\label{e:rescaled-coefficients}
 P_{i,Q}(y)=\sum_{j=1}^{d_i}a_{i,j}Q^{j-1}y^j.
\end{equation}
%The leading coefficient is $a_{i,d_i}Q^{d_i-1}$, and
%$|a_{i,d_i}|\geq1$, so the same $C$ works for every modulus.

The next proposition converts a discrepancy in ordinary-prime counts
into a density increment on an arithmetic progression. The proof closely follows \cite[Section~5]{KMTT}, combining
Proposition~5.2 with the density increment argument in the
proof of Lemma~5.4.

\begin{proposition}
\label{p:KMTT}
Fix $C,d$ and $\ell$. There are constants
$C_1,K_1\geq2$, $0<c_0<1/2$, and $c_1,c_2>0$ with the following
property. Suppose
$N$ is sufficiently large and
\[
 (\log N)^{-1/C_1^2}<\delta\leq1,\qquad
 \left\lceil K_1\delta^{-C_1}\right\rceil\leq s
 \leq(\log N)^{c_0},
\]
and
\[
 Q\leq\exp((\log N)^{2c_0}),\qquad
 P^+(Q)\leq s,\qquad \mathcal L(s)\mid Q.
\]
Let $M=\lfloor(N/Q^{d-1})^{1/d}\rfloor$, and let
$R_0=0,R_1,\ldots,R_\ell$ have distinct increasing degrees,
$(C,Q)$-coefficients, and largest degree $d$. If $A\subset[N]$ has
$|A|=\delta N$ and
\begin{align}\label{e:prime-discrepancy}
 \bigg|&\frac1N\sum_x\E_{y\leq M}\Lambda(Qy+1)
 \prod_{i=0}^{\ell}\1_A(x+R_i(y))\\
 &-\frac1N\sum_x\E_{y\leq M}\Lambda(Qy+1)
 \prod_{i=0}^{\ell}\delta\1_{[N]}(x+R_i(y))\bigg|
 \geq c_1\delta^{\ell+1}\E_{y\leq M}\Lambda(Qy+1),
 \nonumber
\end{align}
then there are integers $b=O_d(1)$, $N'\geq1$, and $a$ such that
\begin{equation}\label{e:density-increment}
 a+Q^b[N']\subset[N],\qquad
 |A\cap(a+Q^b[N'])|\geq(1+c_2)\delta N',
\end{equation}
and, writing $r=\deg R_1$,
\begin{equation}\label{e:increment-length}
 N'\geq Q^{-O_{d,\ell}(1)}\delta^{O_{d,\ell}(1)}
 M^r\exp(-d\sqrt{\log M}).
\end{equation}
\end{proposition}

\begin{proof}
Fix $c_1$ as in \cite[Lemma~5.1]{KMTT}, and let $c_0$ denote
the constant $c'$ chosen in \cite[(5.7)]{KMTT}. Take $C_1,K_1$
as in \cite[(5.8)]{KMTT}, including the enlargements in
\cite[(5.35)]{KMTT}.
By \cite[Proposition~5.2, (5.12)]{KMTT},
hypothesis~\eqref{e:prime-discrepancy} gives an integer $M_0$
satisfying
\begin{align}\label{e:m0}
     M\exp(-\sqrt{\log M})\leq M_0\leq M
\end{align}
for which
\begin{align}\label{e:unweighteddisc}
 \bigg|&\frac1N\sum_{x\in\Z}\E_{y\leq M_0}
 \prod_{i=0}^{\ell}\1_A(x+R_i(y))-\frac1N\sum_{x\in\Z}\E_{y\leq M_0}
 \prod_{i=0}^{\ell}\delta\1_{[N]}(x+R_i(y))
 \bigg|
 \geq \frac{c_1}{4}\delta^{\ell+1}.
\end{align}
Apply \cite[Lemma~5.3]{KMTT} with
$\varepsilon=(c_1/16)\delta^{3(\ell+1)}$ and all inputs equal to $\1_A$. The choices of $C_1,K_1$, together with
$s\geq\lceil K_1\delta^{-C_1}\rceil$ and $\mathcal L(s)\mid Q$,
ensure the divisibility condition \cite[(5.25)]{KMTT}.
% The choices of $C_1,K_1$ ensure the divisibility condition
% \cite[(5.25)]{KMTT}. 
The alternative \cite[(5.26)]{KMTT}
would imply
$\log M_0\ll_{C,d,\ell}(\log N)^{2c_0}+\log\log N=o(\log N)$,
contrary to $\log M_0\geq(3d)^{-1}\log N$ for sufficiently
large $N$.

Thus we obtain partitions $\mathcal B_0,\ldots,\mathcal B_\ell$
of $[N]$ into arithmetic progressions of common difference
$Q^b$, where $1\leq b\ll_d1$. By \cite[(5.28)]{KMTT},
every atom has length at least
\begin{equation}\label{e:resbound}
 |B|\geq Q^{-O_{d,\ell}(1)}
 \delta^{O_{d,\ell}(1)}M_0^{\deg R_1}.
\end{equation}
Let $F_i$ equal the density of $A$ on each atom of
$\mathcal B_i$, extended by zero outside $[N]$.
The approximation \cite[(5.27)]{KMTT} and \eqref{e:unweighteddisc} give
\begin{align}\label{e:localfactorforms}
 \left|\frac1N\sum_x\E_{y\leq M_0}
 \left(
 \prod_{i=0}^{\ell}F_i(x+R_i(y))
 -
 \prod_{i=0}^{\ell}\delta\1_{[N]}(x+R_i(y))
 \right)\right|
 \geq\frac{c_1}{8}\delta^{\ell+1}.
\end{align}
We argue as in \cite[proof of Lemma~5.4, immediately after
(5.43)]{KMTT}. Fix $c_2>0$ so that
$2(\ell+1)c_2(1+c_2)^\ell<c_1/8$, and suppose that
$F_i\leq(1+c_2)\delta$ for every $i$.
Since $\sum_xF_i(x)=\delta N$, the positive and negative
parts of $F_i-\delta\1_{[N]}$ have equal mass, giving
$$\sum_x|F_i(x)-\delta\1_{[N]}(x)|\leq2c_2\delta N.$$
Expanding the difference of products one factor at a time
bounds the left-hand side of~\eqref{e:localfactorforms} by
\[
 \frac{((1+c_2)\delta)^\ell}N
 \sum_{i=0}^{\ell}\sum_x
 |F_i(x)-\delta\1_{[N]}(x)|
 \leq2(\ell+1)c_2(1+c_2)^\ell\delta^{\ell+1}.
\]
The resulting contradiction shows that some $F_i$ exceeds
$(1+c_2)\delta$ on an atom $a+Q^b[N']\subset[N]$.
By the definition of $F_i$, this gives
\eqref{e:density-increment}.

Writing $r=\deg R_1$, the resolution bound from \eqref{e:resbound} gives
$N'\geq Q^{-O_{d,\ell}(1)}\delta^{O_{d,\ell}(1)}M_0^r$.
Combining this with the lower bound for $M_0$ in \eqref{e:m0} and $r\leq d$
proves~\eqref{e:increment-length}.
\end{proof}

The next lemma verifies the discrepancy hypothesis using the absence of
a Piatetski--Shapiro-prime configuration.

\begin{lemma}
\label{l:PS-discrepancy}
Fix the family~\eqref{e:polynomials}.
Choose $C\geq1$ bounding the absolute values of all coefficients
of $P_1,\ldots,P_\ell$, and fix the corresponding constants
$C_1,K_1,c_0,c_1$ from Proposition~\ref{p:KMTT}.
This choice of $C$ gives a common coefficient bound for every
rescaled family~\eqref{e:rescaled-polynomial}.
There exists $c_*=c_*(\deg P_{\ell})>1$ such that the following holds for every
$1<c<c_*$, with $\gamma=1/c$. Let $N$ be a sufficiently large
positive integer in terms of $c,\mathbf P$, and let
$A\subset[N]$ have density $\delta=|A|/N$ satisfying
$(\log N)^{-1/C_1^2}<\delta\leq1$.
Let $s\geq2$ be an integer satisfying
\[
 \left\lceil K_1\delta^{-C_1}\right\rceil
 \leq s\leq(\log N)^{c_0},
\]
and let $Q\geq1$ be an integer satisfying
\[
 Q\leq\exp((\log N)^{2c_0}),\qquad
 P^+(Q)\leq s,\qquad \mathcal L(s)\mid Q.
\]
Put $M=\lfloor(N/Q^{d-1})^{1/d}\rfloor$.
If $A$ contains no configuration of pairwise distinct elements
\begin{align}\label{e:forbidencond}
     x,\quad x+P_{1,Q}(y),\quad\ldots,\quad x+P_{\ell,Q}(y)
\end{align}
with $x\in\Z$, $y\in\N$, and $Qy+1\in\mathbb P_c$, then
\eqref{e:prime-discrepancy} holds with $R_0=0$ and
$R_i=P_{i,Q}$ for $1\leq i\leq\ell$.
\end{lemma}

\begin{proof}
First we establish a lower bound for the ordinary-prime reference
count, with all inputs equal to $\delta\1_{[N]}$.

Choose a sufficiently small fixed $\tau=\tau(\mathbf P)>0$ and put
$M_\tau=\lfloor\tau M\rfloor$. The coefficient bounds and the
definition of $M$ give
\[
 \max_{1\leq i\leq\ell}|P_{i,Q}(y)|\leq N/10
 \qquad(1\leq y\leq M_\tau).
\]
Hence, for each such $y$, at least $N/2$ integers $x$ satisfy
$x+P_{i,Q}(y)\in[N]$ for every $0\leq i\leq\ell$. Therefore
\begin{align}\label{e:model-lower}
 &\frac1N\sum_x\E_{y\leq M}\Lambda(Qy+1)
 \prod_{i=0}^{\ell}\delta\1_{[N]}(x+P_{i,Q}(y))
 \nonumber\\
 &\qquad\geq
 \frac{\delta^{\ell+1}}{2M}\sum_{y\leq M_\tau}\Lambda(Qy+1)
 \gg_{\mathbf P}
 \delta^{\ell+1}\E_{y\leq M}\Lambda(Qy+1).
\end{align}
The last inequality is the ordinary-prime mass comparison
\cite[(5.2)]{KMTT}.

Taking $\mathbf f=(\1_A,\ldots,\1_A)$, the identity
$D_\gamma=\Lambda_\gamma-\Lambda$ and
Proposition~\ref{p:comparison} give
\begin{align}\label{e:comparisonterms}
     \mathcal T_{N,M,Q}(\Lambda;\mathbf f)
 =
 \mathcal T_{N,M,Q}(\Lambda_\gamma;\mathbf f)
 +O_{c,\mathbf P}(N^{-\xi}).
\end{align}
Since $\Lambda_\gamma$ is supported on $\mathcal S_c$, every nonzero
contribution on the right with $Qy+1$ prime comes from $\mathbb P_c$.
The  hypothesis \eqref{e:forbidencond} eliminates these terms whenever
$0,P_{1,Q}(y),\ldots,P_{\ell,Q}(y)$ are pairwise distinct.
So only collisions and higher prime powers remain.

Each polynomial $P_{i,Q}$, $1\leq i\leq\ell$, and each difference
$P_{i,Q}-P_{j,Q}$, $1\leq i<j\leq\ell$, is nonzero and has
degree at most $d$.
Hence collisions occur for only $O_{d,\ell}(1)$ values of $y$.
Their total contribution is
\begin{align}\label{e:collisionterms}
    \ll_{c,\mathbf P}\frac{(QM)^{1-\gamma}\log(2QM)}M.
\end{align}
The contribution of higher prime powers is bounded, using
the elementary estimate in \cite[(5.5)]{KMTT} and
$\Lambda_\gamma(n)\ll_\gamma(QM)^{1-\gamma}\Lambda(n)$
for $n\leq QM+1$, by
\begin{equation}\label{e:prime-powers}
 \ll_\gamma Q^{3/2-\gamma}M^{1/2-\gamma}
 (\log(2QM))^2.
\end{equation}
% The remaining contribution comes from $Qy+1=p^j$ with $j\geq2$. Since
% \[
%  \sum_{\substack{n\leq L+1\\n=p^j,\ j\geq2}}\Lambda(n)
%  \ll L^{1/2}(\log(2L))^2,
% \]
% it is at most
% \begin{equation}\label{e:prime-powers}
%  \frac{(QM)^{1-\gamma}}M
%  \sum_{\substack{n\leq QM+1\\n=p^j,\ j\geq2}}\Lambda(n)
%  \ll Q^{3/2- \gamma}M^{1/2-\gamma}(\log(2QM))^2.
% \end{equation}
% Because $Q=M^{o(1)}$ and $\gamma$ is close to $1$, \eqref{e:collisionterms},
% \eqref{e:prime-powers}, and the error term in \eqref{e:comparisonterms} are all $O(N^{-\eta})$
% for some $\eta>0$.
Because $Q=M^{o(1)}$ and $\gamma$ is close to $1$,
\eqref{e:collisionterms}, \eqref{e:prime-powers}, and the error
term in \eqref{e:comparisonterms} are all $O(N^{-\eta})$
for some $\eta>0$. Finally, the prime-mass estimate\footnote{The recent exclusion of Landau--Siegel zeros,
originating in an OpenAI argument and presented by Green
\cite{GreenNoSiegel} (see also Tao's exposition
\cite{TaoNoSiegel}), combined with the prime number theorem
in arithmetic progressions, gives the stronger bound
$\E_{y\leq M}\Lambda(Qy+1)\gg Q/\varphi(Q)$
in the present parameter range. For the sake of completeness,
we retain the classical bounds in the proof.} \cite[(5.15)]{KMTT} gives
\begin{equation}\label{e:prime-mass}
 \E_{y\leq M}\Lambda(Qy+1)
 \gg\exp\bigl(-O((\log N)^{2c_0})\bigr).
\end{equation}
Thus $N^{-\eta}$ is negligible compared with the right-hand side of
\eqref{e:model-lower}, uniformly in the density range of
Proposition~\ref{p:KMTT}. The constant $c_1$ in
\cite[(5.11)]{KMTT} is chosen as a fixed sufficiently small fraction of
the implicit constant in~\eqref{e:model-lower}. The preceding estimates show that
$\mathcal T_{N,M,Q}(\Lambda;\mathbf f)
=O_{c,\mathbf P}(N^{-\eta})$.
Subtracting this count in $A$ from the  count
in~\eqref{e:model-lower}, we obtain, for sufficiently large $N$,
\begin{align*}
 &\frac1N\sum_x\E_{y\leq M}\Lambda(Qy+1)
 \prod_{i=0}^{\ell}\delta\1_{[N]}(x+P_{i,Q}(y))\\
 &\quad-\frac1N\sum_x\E_{y\leq M}\Lambda(Qy+1)
 \prod_{i=0}^{\ell}\1_A(x+P_{i,Q}(y))\\
 &\qquad\geq
 c_1\delta^{\ell+1}\E_{y\leq M}\Lambda(Qy+1).
\end{align*}
This proves~\eqref{e:prime-discrepancy}.
\end{proof}

\begin{proof}[Proof of Theorem~\ref{t:szemeredi}]
Let $\alpha=|A|/N$ and suppose, towards a contradiction, that $A$ contains
no configuration~\eqref{e:pattern} with pairwise distinct terms.  Let
$C_1,K_1,c_0,c_2$ be supplied by Proposition~\ref{p:KMTT}, and put
\[
 s=\left\lceil K_1\alpha^{-C_1}\right\rceil,
 \qquad Q_0=\mathcal L(s).
\]
Since $Q_0=\mathcal L(s)\leq s!$, we have
\[
 \log Q_0\leq s\log s
 \ll_{\mathbf P}\alpha^{-O_{\mathbf P}(1)}.
\]
The assumption $\alpha\geq(\log N)^{-\kappa}$ therefore gives
$\log Q_0=o(\log N)$, provided $\kappa>0$ is sufficiently small.
Thus $Q_0=N^{o(1)}$.

Choose a residue class modulo $Q_0$ on whose intersection with
$[N]$ the set $A$ has maximal relative density, and represent
it by $r_0\in\{1-Q_0,\ldots,0\}$. Put
\[
 N_0:=\left\lfloor\frac{N-r_0}{Q_0}\right\rfloor,\qquad
 A_0:=\{u\in[N_0]:r_0+Q_0u\in A\}.
\]
Then $r_0+Q_0[N_0]=(r_0+Q_0\Z)\cap[N]$.
Averaging over the full intersections of the residue classes
with $[N]$ gives $\delta_0:=|A_0|/N_0\geq\alpha$.
Moreover, for sufficiently large $N$,
\begin{equation}\label{e:initial-scale}
 N_0\geq N/Q_0-1\geq N/(2Q_0).
\end{equation}

Since $r_0+Q_0A_0\subseteq A$, the forbidden-configuration
hypothesis \eqref{e:pattern} on $A$ implies that $A_0$ contains no configuration
of pairwise distinct elements
\[
 x,\ x+P_{1,Q_0}(y),\ldots,x+P_{\ell,Q_0}(y),
 \qquad y\in\N,\quad Q_0y+1\in\mathbb P_c.
\]
Indeed, if such a configuration existed, setting
$a=r_0+Q_0x$ and $p=Q_0y+1$ would give $a\in A$ and
\[
 a+P_i(p-1)
 =r_0+Q_0\bigl(x+P_{i,Q_0}(y)\bigr)\in A
 \qquad(1\leq i\leq\ell).
\]
Multiplication by $Q_0$ preserves distinctness, so these elements
would form a configuration~\eqref{e:pattern} in $A$, a contradiction.

We iterate Lemma~\ref{l:PS-discrepancy} and
Proposition~\ref{p:KMTT}. Suppose that $A_j\subset[N_j]$ has density
$\delta_j\geq\alpha$ and contains no configuration of pairwise
distinct elements
\[
 x,\ x+P_{1,Q_j}(y),\ldots,x+P_{\ell,Q_j}(y),
 \qquad y\in\N,\quad Q_jy+1\in\mathbb P_c.
\]
The density increment gives a progression
$a_j+Q_j^{b_j}[N_{j+1}]\subset[N_j]$, with $b_j=O_d(1)$, on which
$A_j$ has density at least $(1+c_2)\delta_j$. Put
\[
 A_{j+1}:=\{u\in[N_{j+1}]:a_j+Q_j^{b_j}u\in A_j\},
 \qquad
 \delta_{j+1}:=\frac{|A_{j+1}|}{N_{j+1}},
\]
and set
\begin{equation}\label{e:next-modulus}
 Q_{j+1}:=Q_j^{b_j+1}.
\end{equation}
In particular, the first step produces $A_1$, and
$\delta_{j+1}\geq(1+c_2)\delta_j$ at every step.

The forbidden-configuration \eqref{e:pattern} hypothesis is preserved, since
\[
 Q_j^{b_j}P_{i,Q_{j+1}}(y)
 =P_{i,Q_j}(Q_j^{b_j}y),
 \qquad
 Q_{j+1}y+1=Q_j(Q_j^{b_j}y)+1.
\]
Thus multiplication by $Q_j^{b_j}$ and translation by $a_j$ would
take a forbidden configuration in $A_{j+1}$ to one in $A_j$.

We keep the smoothness parameter $s$ fixed throughout. This is
permitted by Proposition~\ref{p:KMTT}, because
\[
 \left\lceil K_1\delta_j^{-C_1}\right\rceil
 \leq \left\lceil K_1\alpha^{-C_1}\right\rceil=s.
\]
Moreover, every $Q_j$ is a power of $Q_0$, so
\[
 P^+(Q_j)=P^+(Q_0)\leq s\ll_{\mathbf P}\alpha^{-C_1},
 \qquad \mathcal L(s)\mid Q_j.
\]
Hence the smoothness bound remains in terms of the initial density
$\alpha$, even as $\delta_j$ increases.

We now verify inductively that the hypotheses required for each
density increment remain valid.  Since $b_j=O_d(1)$ and there are at most
$O_{\mathbf P}(\log\alpha^{-1}+1)$ possible increments,
\cite[(5.57)]{KMTT} gives
\begin{equation}\label{e:modulus-control}
 \log Q_j\leq\alpha^{-O_{\mathbf P}(1)}.
\end{equation}
Writing $r=\deg P_1$ and $\vartheta=r/d$, the length estimate
\eqref{e:increment-length} gives
\begin{equation}\label{e:length-recursion}
 \log N_{j+1}\geq\vartheta\log N_j
 -O_{\mathbf P}\bigl(\log Q_j+\log(\alpha^{-1})
 +\sqrt{\log N_j}\bigr),
\end{equation}
which is \cite[(5.58)]{KMTT}.  Iterating this inequality, using
\eqref{e:initial-scale}, \eqref{e:modulus-control}, and
$j=O_{\mathbf P}(\log\alpha^{-1}+1)$, gives
\begin{equation}\label{e:length-control}
 \log N_j\gg_{\mathbf P}
 \alpha^{O_{\mathbf P}(1)}\log N
 -\alpha^{-O_{\mathbf P}(1)}
 -O_{\mathbf P}\bigl(\log(\alpha^{-1})\sqrt{\log N}\bigr);
\end{equation}
see \cite[(5.59)]{KMTT}.
If $\alpha\geq(\log N)^{-\kappa}$ and $\kappa>0$ is sufficiently small,
the first term in~\eqref{e:length-control} dominates the errors.  Hence
$\log N_j\geq(\log N)^{\eta_0}$ for some $\eta_0>0$, uniformly at every
stage.  Equations~\eqref{e:modulus-control} and
\eqref{e:length-control} then imply
\[
 Q_j\leq\exp((\log N_j)^{2c_0}),\qquad
 s\leq(\log N_j)^{c_0},\qquad
 \delta_j>(\log N_j)^{-1/C_1^2},
\]
after decreasing $\kappa$ and $c_0$ if necessary.  Since every
$Q_j$ is a power of $Q_0$, one also has
$P^+(Q_j)\leq s$ and $\mathcal L(s)\mid Q_j$.  These are exactly the
hypotheses of Proposition~\ref{p:KMTT}.  Applying
Proposition~\ref{p:comparison} with any fixed $a$ satisfying
$2c_0<a<1$ also puts every $Q_j$ in its required comparison range.
We may therefore perform
\[
 \left\lfloor\frac{\log(\alpha^{-1})}{\log(1+c_2)}\right\rfloor+1
\]
increments.  The resulting density is greater than $1$, a contradiction.
This proves the theorem.
\end{proof}

\section{Pointwise ergodic averages}

By homogeneity, assume throughout this section that
$\|f_i\|_\infty\leq1$. Write
$z_n(x):=\prod_{i=1}^{\ell}f_i(T^{P_i(n-1)}x)$ and put
\begin{equation}\label{e:weighted-averages}
 B_N^\gamma(x):=\E_{n\leq N}\Lambda_\gamma(n)z_n(x),
 \qquad
 B_N(x):=\E_{n\leq N}\Lambda(n)z_n(x).
\end{equation}
Taking the zero phase in Theorem~\ref{t:shao-sun} and using the
prime number theorem gives
\begin{equation}\label{e:weight-mass}
 \E_{n\leq N}\Lambda_\gamma(n)
 =\E_{n\leq N}\Lambda(n)+O_{c,d}(N^{-\sigma_d})
 =1+o(1).
\end{equation}

\begin{proposition}[Pointwise comparison]\label{p:pointwise-comparison}
After decreasing $c_*>1$ if necessary, for every $1<c<c_*$ one has
$B_N^\gamma-B_N\to0$ almost everywhere.
\end{proposition}

\begin{proof}
By definition,
$B_N^\gamma-B_N=\E_{n\leq N}D_\gamma(n)z_n$.
We apply \cite[Proposition~1.4]{mousavisolo} with $g=D_\gamma$
and the polynomials $n\mapsto P_i(n-1)$, which retain their
distinct degrees. Thus, for any $B>0$ satisfying $|D_\gamma(n)|\leq B$ on $[N]$,
\[
 \|B_N^\gamma-B_N\|_{L^1(X)}
 \ll_{\mathbf P}
 B\left(
 N^{-1}+B^{-1}
 \sup_{\deg\phi\leq d}
 \left|\E_{n\leq N}D_\gamma(n)\e(\phi(n))\right|
 \right)^{1/K}.
\]
We now take $B=C_cN^{1-\gamma}\log(2N)$, with $C_c$ sufficiently
large.  Theorem~\ref{t:shao-sun} gives
$$\sup_{\deg\phi\leq d}|\E_{n\leq N}D_\gamma(n)\e(\phi(n))|
\ll_{c,d}N^{-\sigma_d}.$$ 
% while
% $|D_\gamma(n)|\ll_c N^{1-\gamma}\log(2N)$ for $n\leq N$.
Consequently,
\begin{align}
 \|B_N^\gamma-B_N\|_{L^1(X)}
 &\ll_{c,\mathbf P}
 N^{1-\gamma}\log(2N)
 \left(N^{-1}+N^{\gamma-1-\sigma_d}\right)^{1/K}
 \nonumber\\
 &\ll_{c,\mathbf P}\log(2N)
 \left(
 N^{1-\gamma-1/K}
 +N^{(1-\gamma)(1-1/K)-\sigma_d/K}
 \right)
\ll_{c,\mathbf P}N^{-\eta}.
 \label{e:L1-comparison}
\end{align}
Here, after decreasing $c_*$ so that
$1-\gamma\leq\min\{1,\sigma_d\}/(2K)$, we may take
$\eta=\min\{1,\sigma_d\}/(4K)$.

% We use the orbit estimate \cite[Proposition~1.4]{mousavisolo},
% applied to the polynomials $n\mapsto P_i(n-1)$, which retain
% their distinct degrees. Together with Theorem~\ref{t:shao-sun}
% and the bound
% $|D_\gamma(n)|\ll_c N^{1-\gamma}\log(2N)$ for $n\leq N$, it gives
% \begin{align}
%  \|B_N^\gamma-B_N\|_{L^1(X)}
%  &\ll_{c,\mathbf P}
%  N^{1-\gamma}\log(2N)
%  \left(N^{-1}+N^{\gamma-1-\sigma_d}\right)^{1/K}
% \ll_{c,\mathbf P}N^{-\eta},
%  \label{e:L1-comparison}
% \end{align}
% for some $\eta>0$, provided $c_*$ is sufficiently close to $1$.

Choose $A>\max\{1,\eta^{-1}\}$ and put $N_j=\lfloor j^A\rfloor$.
By~\eqref{e:L1-comparison},
$$\sum_j\|B_{N_j}^\gamma-B_{N_j}\|_1<\infty.$$
Tonelli's theorem therefore gives
$B_{N_j}^\gamma(x)-B_{N_j}(x)\to0$ for almost every $x$.

To pass to all $N$, let $w$ be either $\Lambda_\gamma$ or $\Lambda$.
Since $|z_n|\leq1$, for $N_j\leq N<N_{j+1}$ the corresponding
weighted average differs from its value at $N_j$ by at most
\[
 \frac{N_{j+1}-N_j}{N_j^2}\sum_{n\leq N_j}w(n)
 +\frac1{N_j}\sum_{N_j<n\leq N_{j+1}}w(n)
 =o(1),
\]
uniformly in $x$. Here we used~\eqref{e:weight-mass} and
$N_{j+1}/N_j\to1$. Applying this to both weighted averages
proves the full-sequence comparison.
\end{proof}

\begin{proof}[Proof of Theorem~\ref{t:pointwise}]
By \cite[Theorem~1.1(ii)]{MirekWanWright}, applied to the
polynomials $n\mapsto P_i(n-1)$, the averages $B_N$ converge
almost everywhere. Proposition~\ref{p:pointwise-comparison}
gives the same limit, denoted by $L(x)$, for $B_N^\gamma$.

We remove the weights by summation by parts. First, higher
prime powers contribute at most
\begin{align}\label{e:primepowers}
    \sum_{\substack{p^k\leq N\\k\geq2}}
 \bigl(\Lambda_\gamma(p^k)+\Lambda(p^k)\bigr)
 \ll_\gamma N^{3/2-\gamma}(\log(2N))^2=o(N),
\end{align}
since $\gamma>1/2$ after decreasing $c_*$ if necessary.

Fix $x$ for which $B_N^\gamma(x)\to L(x)$, and put
$a(t)=\gamma^{-1}t^{1-\gamma}\log t$ for $t\geq2$. For a prime $p$,
we have $\Lambda_\gamma(p)=a(p)\1_{\mathbb P_c}(p)$. Thus removing
higher prime powers leaves precisely the sums over $\mathbb P_c$.

Define $W(t):=\sum_{p\leq t,\ p\in\mathbb P_c}a(p)$.
Estimate \eqref{e:primepowers}, together with $|z_n(x)|\leq1$,
gives
\begin{align*}
 W(N)
 &=\sum_{n\leq N}\Lambda_\gamma(n)+o(N)
 =N+o(N),\\
 \sum_{\substack{p\leq N\\p\in\mathbb P_c}}a(p)z_p(x)
 &=N B_N^\gamma(x)+o(N)
 =N L(x)+o_x(N).
\end{align*}
The first asymptotic uses~\eqref{e:weight-mass}; the second uses
the convergence of $B_N^\gamma(x)$. Subtracting $L(x)W(N)$
from the second sum therefore yields
\[
 \sum_{\substack{p\leq N\\p\in\mathbb P_c}}
 a(p)\bigl(z_p(x)-L(x)\bigr)
 =N L(x)+o_x(N)-L(x)\bigl(N+o(N)\bigr)
 =o_x(N).
\]
Since $W(N)\sim N$, the centered weighted sum is also
$o_x(W(N))$. Write
\[
 S(t):=\sum_{\substack{p\leq t\\p\in\mathbb P_c}}
 a(p)\bigl(z_p(x)-L(x)\bigr).
\]
For real $t\geq2$, the conditions $p\leq t$ and
$p\leq\lfloor t\rfloor$ are equivalent. Thus
$S(t)=S(\lfloor t\rfloor)$ and $W(t)=W(\lfloor t\rfloor)$.
Consequently, given $\varepsilon>0$, there is an integer
$t_0\geq2$ such that $|S(t)|\leq\varepsilon W(t)$ for every
real $t\geq t_0$.
For $N\geq t_0$, summation by parts and $a'(t)\geq0$ give
\begin{align*}
 \left|\sum_{\substack{p\leq N\\p\in\mathbb P_c}}
 \bigl(z_p(x)-L(x)\bigr)\right|
 &=\left|
 \frac{S(N)}{a(N)}
 +\int_2^N S(t)\frac{a'(t)}{a(t)^2}\,dt
 \right|\\
 &\leq C_{x,\varepsilon}
 +\varepsilon\left(
 \frac{W(N)}{a(N)}
 +\int_2^N W(t)\frac{a'(t)}{a(t)^2}\,dt
 \right)\\
 &=C_{x,\varepsilon}
 +\varepsilon\sum_{\substack{p\leq N\\p\in\mathbb P_c}}
 \frac{a(p)}{a(p)}\\
 &=C_{x,\varepsilon}+\varepsilon\pi_c(N).
\end{align*}
Here
$C_{x,\varepsilon}:=\int_2^{t_0}|S(t)|a'(t)/a(t)^2\,dt$
bounds the initial part of the integral. The equality involving
$W$ follows by applying the same summation-by-parts formula
to its defining sum.

Since $W(N)\sim N$, the set $\mathbb P_c$ is infinite and
$\pi_c(N)\to\infty$. Dividing by $\pi_c(N)$ and then letting
$\varepsilon\to0$ proves
$A_N^{\mathbb P_c}(\mathbf f)(x)\to L(x)$.

Applying the same summation-by-parts calculation over all primes,
with $a(t)=\log t$ and
$W(t)=\sum_{p\leq t}\log p\sim t$, shows that
the ordinary-prime average~\eqref{e:prime-average} also converges
to $L(x)$. This proves the theorem.
\end{proof}

For $1\leq r<\infty$ and a sequence $(a_j)_{j\geq1}$ of
complex numbers, define its $r$-variation seminorm by
\[
 V^r(a_j:j\geq1)
 :=\sup_{\substack{J\geq1\\1\leq j_0<\cdots<j_J}}
 \left(\sum_{h=1}^{J}
 |a_{j_h}-a_{j_{h-1}}|^r\right)^{1/r}.
\]
For a sequence of functions, this definition is applied pointwise.
\begin{corollary}
\label{c:long-variation}
Under the hypotheses of Theorem~\ref{t:pointwise}, let $r>2$
and $\lambda>1$. If $(N_j)_{j\geq1}$ is a sequence of positive
integers satisfying $N_{j+1}\geq\lambda N_j$, then
\[
 \left\|V^r\bigl(B_{N_j}^\gamma:j\geq1\bigr)\right\|_{L^1(X)}
 \ll_{c,\mathbf P,r,\lambda}
 \prod_{i=1}^{\ell}\|f_i\|_\infty.
\]
\end{corollary}

\begin{proof}
By homogeneity, assume that $\|f_i\|_\infty\leq1$.
The elementary bound $V^r(a_j:j\geq1)\leq2\sum_j|a_j|$,
together with~\eqref{e:L1-comparison}, gives
\begin{align*}
 \left\|V^r\bigl(B_{N_j}^\gamma-B_{N_j}:j\geq1\bigr)\right\|_1
 &\leq2\sum_{j\geq1}
 \|B_{N_j}^\gamma-B_{N_j}\|_1\ll_{c,\mathbf P}\sum_{j\geq1}N_j^{-\eta}
 \ll_{c,\mathbf P,\lambda}1.
\end{align*}
Since $\nu(X)=1$, \cite[Theorem~1.1(iv)]{MirekWanWright}
and H\"older's inequality give
$$\|V^r(B_{N_j}:j\geq1)\|_1\ll_{\mathbf P,r,\lambda}1.$$
The triangle inequality completes the proof.
\end{proof}


\begin{thebibliography}{99}
\raggedright

\bibitem{AkbalGuloglu}
Y.~Akbal and A.~M.~G\"ulo\u{g}lu,
Waring--Goldbach problem with Piatetski--Shapiro primes,
\emph{J. Th\'eor. Nombres Bordeaux} \textbf{30} (2018),
no.~2, 449--467.


\bibitem{BalogFriedlander}
A.~Balog and J.~B.~Friedlander,
A hybrid of theorems of Vinogradov and Piatetski--Shapiro,
\emph{Pacific J. Math.} \textbf{156} (1992), no.~1, 45--62.

\bibitem{BergelsonLeibman}
V. Bergelson and A. Leibman,
\emph{Polynomial extensions of van der Waerden's and Szemer\'edi's
theorems},
J. Amer. Math. Soc. \textbf{9} (1996), no.~3, 725--753.

\bibitem{BourgainPrimes}
J. Bourgain,
\emph{An approach to pointwise ergodic theorems},
in Geometric Aspects of Functional Analysis (1986/87),
Lecture Notes in Math. \textbf{1317},
Springer, Berlin, 1988, 204--223.

\bibitem{BurginDigits}
A. Burgin,
\emph{S\'ark\"ozy's theorem for shifted primes with restricted digits},
arXiv:2510.13076 (2025).

\bibitem{FornalKrause}
J. Fornal and B. Krause,
\emph{Double recurrence and almost sure convergence:
primes and weighted theory},
arXiv:2603.22203.

\bibitem{GreenNoSiegel}
B.~Green,
\emph{On Landau--Siegel zeros},
Hexagon preprint 2610.00311v2 (2026),
\url{https://hexagonmath.org/2610.00311v2}.

\bibitem{KMTTPointwise}
B. Krause, H. Mousavi, T. Tao and J. Ter\"av\"ainen,
\emph{Pointwise convergence of bilinear polynomial averages over the primes},
Ergodic Theory Dynam. Systems (2025),
doi:10.1017/etds.2025.10202; arXiv:2409.10510.


\bibitem{KMTT}
B. Krause, H. Mousavi, T. Tao and J. Ter\"av\"ainen,
\emph{Quantitative bounds for sets lacking polynomial progressions with
shifted prime difference},
arXiv:2608.19525v1 (2026).

\bibitem{LaVictoireBad}
P.~LaVictoire,
\emph{Universally $L^1$-bad arithmetic sequences},
J. Anal. Math. \textbf{113} (2011), 241--263.

\bibitem{LiPanPS}
H.~Li and H.~Pan,
\emph{The Green--Tao theorem for Piatetski--Shapiro primes},
J. Funct. Anal. \textbf{291} (2026), no.~1,
Article~111440, doi:10.1016/j.jfa.2026.111440.


\bibitem{MaynardDigits}
J. Maynard,
\emph{Primes with restricted digits},
Invent. Math. \textbf{217} (2019), no.~1, 127--218.


\bibitem{MirekPS}
M. Mirek,
\emph{Roth's theorem in the Piatetski--Shapiro primes},
Rev. Mat. Iberoam. \textbf{31} (2015), no.~2, 617--656.

\bibitem{MirekThinPrimes}
M. Mirek,
\emph{$\ell^p(\mathbb Z)$-boundedness of discrete maximal functions
along thin subsets of primes and pointwise ergodic theorems},
arXiv:1401.1361.

\bibitem{MirekWanWright}
M. Mirek, R. Wan and J. Wright,
\emph{Pointwise convergence of polynomial multiple ergodic averages along
the primes},
arXiv:2505.15549v2 (2026).



\bibitem{mousavisolo} H. Mousavi, \emph{Polynomial configurations and pointwise averages along Piatetski-Shapiro sequences}, arXiv:2609.01847. 2026 Sep 1.

\bibitem{MousaviPSFormalization}
H.~Mousavi,
\emph{Lean formalization of polynomial configurations and pointwise
averages along Piatetski--Shapiro primes},
\href{https://github.com/HamedMousaviMath/Lean-Formalization-Polynomial-Configurations-and-Pointwise-Averages-Along-Piatetski--Shapiro-Primes}{GitHub repository link}.

\bibitem{PiatetskiShapiro}
I.~I.~Piatetski--Shapiro,
\emph{On the distribution of prime numbers in sequences of
the form $[f(n)]$},
Mat. Sb. (N.S.) \textbf{33(75)} (1953), no.~3, 559--566
(in Russian).




\bibitem{SarkozySquares}
A.~S\'ark\"ozy,
\emph{On difference sets of sequences of integers. I},
Acta Math. Acad. Sci. Hungar. \textbf{31} (1978),
no.~1--2, 125--149.

\bibitem{ShaoWang}
X.~Shao and M.~Wang,
\emph{Quantitative bounds in a popular polynomial Szemer\'edi theorem},
Proc. Roy. Soc. Edinburgh Sect. A (2025), First View, 1--27,
doi:10.1017/prm.2025.10103.

\bibitem{ShaoSun}
X. Shao and Y.-C. Sun,
\emph{Linear equations in Piatetski--Shapiro primes},
arXiv:2605.18676v1 (2026).


\bibitem{TaoNoSiegel}
T.~Tao,
\emph{A sketch of the proof of non-existence of Siegel zeroes},
What's new, 19 October 2026,
\url{https://terrytao.wordpress.com/2026/10/19/a-sketch-of-the-proof-of-non-existence-of-siegel-zeroes/}.

\bibitem{Teravainen}
J. Ter\"av\"ainen,
\emph{Pointwise convergence of ergodic averages with M\"obius weight},
arXiv:2401.03174v2 (2024).

\bibitem{WierdlPrimes}
M. Wierdl,
\emph{Pointwise ergodic theorem along the prime numbers},
Israel J. Math. \textbf{64} (1988), 315--336.

\bibitem{WooleyZiegler}
T. D. Wooley and T. D. Ziegler,
\emph{Multiple recurrence and convergence along the primes},
Amer. J. Math. \textbf{134} (2012), no.~6, 1705--1732.






\end{thebibliography}
\end{document}